\documentclass[12pt,reqno]{article}
\usepackage{amsfonts}
\usepackage{graphicx}
\usepackage{epstopdf}
\usepackage{amsmath}
\usepackage{amsthm}
\usepackage{latexsym,bm}
\usepackage{amssymb}
\usepackage{indentfirst}
\newtheorem{thm}{Theorem}[section]
\newtheorem{cor}[thm]{Corollary}
\newtheorem{lem}[thm]{Lemma}

\newtheorem{defn}[thm]{Definition}
\newtheorem{rem}[thm]{Remark}
\newtheorem{exm}{Example}

\numberwithin{equation}{section}

\begin{document}

\title{\bf A  conformal characterization  of  manifolds of constant sectional curvature}

\author{\bf Xiaoyang Chen\thanks{Partially supported by National Natural Science Foundation of China No.11701427 and the Fundamental Research Funds for the Central Universities.} ,
Francisco Fontenele and Frederico Xavier\thanks{Partially supported by the John William and Helen Stubbs Potter Professorship in Mathematics.}}

\date{}
\maketitle

%

\begin{quote}
\small {\bf Abstract}.
A special case of the main result states that  a  complete  $1$-connected Riemannian manifold $(M^n,g)$  is isometric to  one of the  models  $\mathbb R^n$, $S^n(c)$,
$\mathbb H^n(-c)$  of constant curvature if and only if  every $p\in M^n$  is a non-degenerate maximum of  a      germ of smooth functions whose Riemannian gradient is a conformal vector field.

\end{quote}

\section{Introduction}

The goal of  the present  paper is to  offer a  characterization of   connected Riemannian manifolds of constant sectional curvature in terms of the existence of certain special germs of functions.
For conceptual clarity,   it is convenient to single out the following notion:
\begin{defn} \rm
A $\underline{ \rm conformal \; Morse \; germ} $    (CMG, for short) based at a point $p$ of a Riemannian manifold $(M^n,g)$ is a germ $[f]_p$ of  smooth   functions  for which $p$ is a non-degenerate critical point of (any representative) $f$ and, relative to  $g$,     the gradient $\nabla f$ is a conformal vector field on a neighborhood of $p$. \end{defn}

Recall that conformality of $\nabla f$  has a dynamical meaning:  for all small $|t|$,  the time-$t$ map of the local flow of
$\nabla f$ is conformal in the usual sense, i.e. it is angle-preserving  relative to the metric $g$.
It is easy to see that $[f]_p$ is a CMG in $(M^n,g)$ if and only if $\nabla f(p)=0$ and $\nabla^2f=hg$ for some function $h$ with $h(p)\neq 0$. It follows that  if $p$ is the base of a CMG, then $p$ is either a
 non-degenerate  local maximum or  a non-degenerate local minimum of $f$.

If $[f]_p$ is a CMG,  so is $[-f]_p$. In particular,  we may assume that $p$ is a local maximum   of $f$, in which case the  forward   flow of $\nabla f$  pushes   all nearby points towards the point of maximum $p$,   while preserving angles between curves.
The existence of this special   dynamics suggests that, at the point $p$ itself, the possibilities for the geometry of $(M^n, g)$ should be severely restricted. The theorem below shows that this is indeed the case:
\begin{thm} \label {point} Let $p$ be the base of a conformal Morse germ in $(M^n,g)$.
\vskip3pt
\noindent i) If $n=2$, then $p$ is a critical point of the curvature.
\vskip3pt
\noindent ii) If $n>2$, then the sectional curvatures of all $2$-planes contained in $T_pM^n$ are equal.
\end{thm}

\noindent \begin{rem}  \rm (Pointwise constancy of sectional curvatures.)  To the best of our knowledge, Theorem \ref{point} ii) is  the only known   criterion for constancy of the  sectional curvatures of all $2$-planes contained in  a $\underline{\rm fixed}$,  preassigned, tangent space $T_pM^n$. In dimensions three and higher, if  this ``pointwise" constancy of the sectional curvatures holds everywhere on the manifold, then $(M^n, g)$  actually has constant curvature (Schur's theorem \cite{KN}), but situations may arise in which one would like to establish that the sectional curvatures are the same at only  a single  point. \end{rem}

\begin{exm} \rm \label{exemplo} One can construct conformal Morse germs on any of the  model spaces $\mathbb R^n$, $S^n(c)$, $\mathbb H^n(-c)$ of constant sectional curvature,
where $c$ is a positive constant. Explicitly, in $S^n(c)$ the function $f(x)=\cos(\sqrt{c}d(x,p))$ satisfies $\nabla^2 f=-cfg$, in $\mathbb R^n$ the function $f(x)=||x-p||^2$ satisfies $\nabla^2 f=2g$, and in $\mathbb H^n(-c)$ the function $f(x)=\cosh(\sqrt{c}d(x,p))$ satisfies $\nabla^2 f=cfg$.
\end{exm}

\noindent The main application of Theorem \ref{point} is the  following  characterization of  manifolds of constant sectional curvature, in terms of conformal Morse germs:
\begin{thm} \label {fundamental criterion}
A connected Riemannian manifold  $(M^n,g)$ has  constant  sectional curvature if and only if  every $p\in M^n$  is a non-degenerate critical point  of a smooth  germ $[f]_p$  for which $\nabla f$ is a conformal vector field.   \end{thm}

\begin{cor} \label{space form}  A  complete  $1$-connected Riemannian manifold $(M^n,g)$  is isometric to  one of the  model spaces  $\mathbb R^n$, $S^n(c)$,
$\mathbb H^n(-c)$ if and only if  every $p\in M^n$  is the base of a CMG.
\end{cor}

Theorem \ref{fundamental criterion} is a direct consequence  of  connectedness, Theorem \ref{point}, Schur's theorem \cite{KN} and, for the converse, Example \ref{exemplo}.
The proof of Theorem \ref{point} itself is rather  involved,   and will be supplied in the next section.

\par The sectional curvatures of  $(M^n,g)$ are not necessarily constant if
 for every point $p \in M^n$  there is a germ $[X]_p$ of conformal vector fields such that $p$ is an isolated non-degenerate zero of  $X$.
 In fact, every Riemannian metric on a two dimensional manifold is locally conformally flat, and so there are always germs  $[X]_p$ as above.  It is only when $X$ is a gradient field that   the curvature is constant (Theorem \ref{fundamental criterion}).

 \par

\begin{rem} \rm In the study of manifolds of constant sectional curvature  it is usually the case that one is forced   to consider  the trichotomy   given by the sign of the curvature (positive, zero, or negative), leading to a separate analysis in each alternative.  For instance, the proof of Obata's theorem (\cite{O}, \cite{BGM}), characterizing the Euclidean sphere, revolves around the equation $\nabla^2f=-fg$, whereas  similar investigations for hyperbolic space would involve the equation $\nabla^2f=fg$.
By contrast,  an appealing feature of Theorem \ref{fundamental criterion} is that it  identifies a relatively simple property that captures the notion of ``constant curvature"  in all three cases simultaneously, independently  of the actual value of the curvature. \end{rem}

We close this Introduction by drawing attention to the  possibility of using the theory of quasiconformal mapings \cite{GMP}  to extend the results in this paper.  Given a point $p$ in a Riemannian manifold $(M, g)$ and $\sigma\subset T_pM$   a $2$-plane, we write $K(p, \sigma)$ for the corresponding sectional curvature. Consider the oscillation of  the  restriction of the sectional curvature function to the Grassmannian
$G(2, T_pM)$:
 \begin{eqnarray} \text{Osc} K_p=   \text{max}_{\sigma\subset T_pM}  K(p, \sigma)-\text{min}_{\sigma\subset T_pM} K(p, \sigma). \end{eqnarray}

 Suppose  that $[f]_p$ is  a germ that has a non-degenerate critical point at $p$   and that, in addition,  $\nabla f$ is $\kappa$-quasiconformal for some $\kappa\geq 1$, in the sense that
  $\kappa-1$ measures the deviation of  $\nabla^2f$ from being a functional multiple of the metric tensor.  Can one use  $[f]_p$ as a probe, so to speak, in order  to estimate the oscillation $ \text{Osc} K_p$  of the sectional curvatures  in terms of the quasiconformality coefficient $\kappa$?  If so,  Theorem \ref{point} ii)  would be a special case of such a result, as  the former  states that  $\text{Osc} K_p=0$ if  $\kappa=1$.

For other  characterizations of spaces of constant curvature see, for instance,  \cite{BGM}, \cite{O} - \cite {WY} and the references therein.

\section{Proof of Theorem \ref{point}}

\noindent The lemma below will be used in the proofs of both halves of Theorem \ref{point}.

\begin{lem} \label{onlyone}  Let $(M,g)$ be an $n$-dimensional  Riemannian manifold,  $p\in M$, $[f]_p$ a    germ for which $p$ is a non-degenerate critical point, and $v\in T_pM$ a unit vector. Then there exists a sequence $(q_k)$ in $M-\{p\}$ such that
$$\lim_{k \to \infty} q_k=p, \;\;\;
\lim_{k \to \infty}\frac{\nabla f(q_k)}{||\nabla  f(q_k)||} = v.$$
\end{lem}
\vskip10pt

To begin the proof of the lemma, let $k\in \{0, \dots, n\}$  be the Morse index of the singularity $p$ of the vector field $X=\nabla f$. By Morse's lemma \cite{H}, one can choose coordinates $(x_1, \dots, x_n)$ so that $0$ corresponds to $p$ and, locally,
$$f(x_1,\dots,x_n)=-x_1^2-\dots-x_k^2+x_{k+1}^2+\dots+x_n^2.$$
Consider the flat metric $g_0=dx_1^2+\dots+dx_n^2$ on a suitably small neighborhood of $p$, and write
$\nabla f$, $\nabla_0f$ for the gradients of $f$ relative to the metrics $g$ and $g_0$, respectively.

Consider the local matrix approximation
\begin{eqnarray*} (g^{ij}(x))=(g^{ij}(0))+ (E^{ij}(x)), \;\; \lim_{x\to 0}|| E^{ij}(x))||=0, \end{eqnarray*}
and let $(g^{ij}(0))=U^{-1}DU$ be a diagonalization of $(g^{ij}(0))$, where $D=\text{diag}(c_1, \dots, c_n)$, $c_j>0$.
For $t\in [0,1]$, write
\begin{eqnarray*} D_t=\text{diag}((1-t)c_1+t, \dots, (1-t)c_n+t), \;\; H(t,x)=U^{-1}D_tU+(1-t)(E^{ij}(x)),\end{eqnarray*}
and observe that the symmetric matrices $H(t,x)$ satisfy
\begin{eqnarray*} H(0,x)=(g^{ij}(x)), \;\;\; H(1,x)=I. \end{eqnarray*}
Since the eigenvalues of $D_t$ are bounded from below on $[0, 1]$ by a positive constant, over a sufficiently small $g_0$-ball $B_{\delta}$ centered at $0$,  one has
\begin{eqnarray*} \inf_{t\in [0,1], x\in B_{\delta}} \text{det} H(t,x) >0. \end{eqnarray*}
In coordinates, the $g$-gradient of $f$ is
\begin{eqnarray*} \nabla f=\sum_{i=1}^n\big (\sum_{j=1}^n g^{ij}(x) \frac {\partial f}{\partial x_j}\big )\frac{\partial}{\partial x_i}.  \end{eqnarray*}
Writing $X$ for the column vector $(\frac {\partial f}{\partial x_1}, \dots, \frac {\partial f}{\partial x_n})^t$, the above expression for $\nabla f$ can be interpreted as the matrix product
\begin{eqnarray*} \nabla f(x)= (g^{ij}(x))X. \end{eqnarray*}
In particular, $H(t,x) X$, with $t\in [0,1]$,
provides a continuous deformation of  the vector field $\nabla f(x)=(g^{ij}(x))X=H(0,x)X$ into
$\nabla_0f=X=H(1, x)X$, through vector fields that have no zeros outside $0$. By the invariance of the Poincar\'e-Hopf index under homotopies \cite{H},
we have
\begin{eqnarray} \label{index} \text{Ind}(\nabla f,p)=\text{Ind}(\nabla_0 f,p)=\text{Ind}((-2x_1,\cdots,-2x_k,2x_{k+1}, \cdots,2x_n),0)
=(-1)^k.\end{eqnarray}
For the convenience of the reader, we briefly recall the definition of Poincar\'e-Hopf index of an isolated singularity.
Let $Y$ be a vector field on $M$ which has an isolated singularity at $p$. Suppose that $\phi: U \rightarrow \phi(U)\subset \mathbb{R}^n$ is a local parametrization
of $M$ carrying $p$ to the origin of $\mathbb{R}^n$, where $U$ is an open set containing $p$. Then the Poincar\'e-Hopf index of $Y$ at $p$ is defined to be
the Poincar\'e-Hopf index of $\widetilde{Y}=d \phi (Y_{|U})$ at the origin, which is equal to the degree of the following map:
$$ F_{\epsilon}: \mathbb{S}_{\epsilon} \rightarrow \mathbb{S}^{n-1}, y \rightarrow \frac{\widetilde{Y}(y)}{\|\widetilde{Y}(y)\|},$$
where $ \mathbb{S}_{\epsilon}  $ is a sphere of radius $\epsilon$ centered at origin. Let $D_0 \widetilde{Y}$ be the differential of
the map $\widetilde{Y}: \phi (U) \rightarrow \mathbb{R}^n$ at the origin. If $D_0 \widetilde{Y}$ is an isomorphism, then
$p$ is said to be a non-degenerate singularity and we have
$$ \text{degree }(F_{\epsilon})=\text{sign}(\text{det} (D_0 \widetilde{Y}) ) \in \{{\pm 1}\}. $$

\par Now,  (\ref{index}) follows from the above description of the Poincar\'e-Hopf index and
the lemma follows from (\ref{index}) also the fact that a continuous non-surjective map $S^{n-1}\to S^{n-1}$ has degree zero.
\qed

\vskip10pt
\noindent To prove Theorem \ref{point}, observe that for all smooth vector fields $X,Z$ one has
\begin{eqnarray}
\nabla^3 f(Z,X)=(\nabla_Z(\nabla^2 f))(X)&=&\nabla_Z\nabla^2 f(X)-\nabla^2 f(\nabla_ZX)\nonumber=\nabla_Z\nabla_X\nabla f-\nabla_{\nabla_ZX}\nabla f.\nonumber
\end{eqnarray}
Likewise,
$ \nabla^3 f(X,Z)=\nabla_X\nabla_Z\nabla f-\nabla_{\nabla_XZ}\nabla f$, and so
\begin{eqnarray}
\nabla^3 f(Z,X)-\nabla^3 f(X,Z)&=&\nabla_Z\nabla_X\nabla f-\nabla_X\nabla_Z\nabla f-\nabla_{[Z,X]}\nabla f\nonumber=R(Z,X)\nabla f,\nonumber
\end{eqnarray}
where $R$ stands for the curvature tensor.
In particular, for any $q\neq p$, $q$ sufficiently close to $p$,  and  unit vector $z$ perpendicular to $\nabla f(q)$, one has
\begin{eqnarray}
\frac{1}{|\nabla f(q)|}\{\nabla^3 f(z,\nabla f(q))-\nabla^3 f(\nabla f(q),z)\}=R(z,\frac{\nabla f(q)}{|\nabla f(q)|})\nabla f(q).\nonumber
\end{eqnarray}
Dividing  by $|\nabla f(q)|$  and  taking the inner product with $z$, one obtains an  expression for the sectional curvature of  the plane spanned by $\nabla f(q)/||\nabla f(q)||$ and $z$:
\begin{eqnarray} \label{Longo}
K(\frac{\nabla f(q)}{|\nabla f(q)|},z)=\frac{1}{|\nabla f(q)|^2}\langle\nabla^3 f(z,\nabla f(q))-\nabla^3 f(\nabla f(q),z),z\rangle.
\end{eqnarray}
Since $\nabla f$ is a conformal vector field, there exists a smooth function $h$ such that $\nabla^2f(Y)=hY$ for every vector field $Y$. Hence, for  $q\neq p$ and any unit vector $z$, $z\perp \nabla f(q)$, one has
\begin{eqnarray}
\nabla^3 f(z,\nabla f(q))=\nabla_z\nabla^2 f(\nabla f)-\nabla^2 f(\nabla_z\nabla f)=\nabla_z(h\nabla f)-h(q)\nabla_z\nabla f=z(h)\nabla f(q).\nonumber
\end{eqnarray}
Similarly,
$\nabla^3 f(\nabla f(q),z)=\nabla f(q)(h)z$.
Substituting in (\ref{Longo}) the values for $\nabla^3 f(z,\nabla f(q))$ and $\nabla^3 f(\nabla f(q),z)$ obtained above, one obtains
\begin{eqnarray} \label{conf}
K(\frac{\nabla f(q)}{|\nabla f(q)|},z)=-\frac{\langle\nabla f,\nabla h\rangle}{|\nabla f|^2}(q).
\end{eqnarray}

For future reference observe that, crucially, the expression for  the sectional curvature of the plane generated by $\nabla f(q)/|\nabla f(q)|$ and $z$  given in (\ref{conf}) is  actually independent of the unit vector   $z$ perpendicular to $\nabla f(q)$.

Assume now that $\text{dim M}=2$ and let $\varphi=\frac{1}{2}|\nabla f|^2$. Since $\nabla^2 f=hI$ for some smooth function $h$, one has
\begin{eqnarray*}
\langle\nabla\varphi,Y\rangle=\frac{1}{2}Y\langle\nabla f,\nabla f\rangle=\langle\nabla_Y\nabla f,\nabla f\rangle=\langle hY,\nabla f\rangle \nonumber
\end{eqnarray*}
for all vector fields $Y$, and so $\nabla\varphi=h\nabla f$. Hence,
\begin{eqnarray*}
\text{Hess}\,\varphi (X,Y)=\langle\nabla_X\nabla\varphi,Y\rangle=\langle X(h)\nabla f+h^2X,Y\rangle=X(h)Y(f)+h^2\langle X,Y\rangle,\nonumber
\end{eqnarray*}
for all  $X,Y$. Similarly,
$\text{Hess}\,\varphi (Y,X)=Y(h)X(f)+h^2\langle Y,X\rangle$.
It follows from the last two equalities, together with  the symmetry of the Hessian,  that
\begin{eqnarray*}
X(h)Y(f)=Y(h)X(f),\;\;\;\forall X,Y,\nonumber
\end{eqnarray*}
and so $X(h)\nabla f=X(f)\nabla h$ for all $X$. Taking $X=\nabla f$, one concludes that
\begin{eqnarray*}
\nabla h(q)=\frac{\langle\nabla h(q),\nabla f(q)\rangle}{|\nabla f(q)|^2}\nabla f(q),\;\;\;q\neq p.\nonumber
\end{eqnarray*}
In view of (\ref{conf}), this yields
\begin{eqnarray}\label{intermediaria}
\nabla h=-K\nabla f.
\end{eqnarray}
Next, we compute the Hessian of $h$. From (\ref{intermediaria}) and $\nabla^2 f=hI$ one obtains
\begin{eqnarray}
\text{Hess}\,h(X,Y)=\langle\nabla_X\nabla h,Y\rangle&=&-\langle\nabla_X(K\nabla f),Y\rangle\nonumber\\&=&-\langle X(K)\nabla f+K\nabla_X\nabla f,Y\rangle\nonumber\\&=&-X(K)Y(f)-Kh\langle X,Y\rangle.\nonumber
\end{eqnarray}
Likewise,
$
\text{Hess}\,h(Y,X)=-Y(K)X(f)-Kh\langle Y,X\rangle$.
Comparing the last two equations  and using the symmetry of the Hessian, one has
\begin{eqnarray}
-X(K)Y(f)=-Y(K)X(f),\;\;\;\forall X,Y,\nonumber
\end{eqnarray}
and  so  $X(K)\nabla f=X(f)\nabla K$ for all vector fields $X$. Taking $X=\nabla f$ in this equality, one concludes that
\begin{eqnarray}\label{final}
\nabla K(q)=\langle\nabla K(q),\frac{\nabla f(q)}{|\nabla f(q)|}\rangle \frac{\nabla f(q)}{|\nabla f(q)|},\;\;\;q\neq p.
\end{eqnarray}
As $\nabla K(q)$ and
$\nabla f(q)$ are collinear by (\ref{final}), and $\nabla f(q)$ approaches  all directions as  $q\to p$ by  Lemma \ref{onlyone}, one must have  $\nabla K(p)=0$. This concludes the proof of Theorem \ref{point} i).
\vskip10pt

We now turn our attention to the second half of Theorem \ref{point}, and assume $\text{dim}\,M>2$. Let $\{e_1,...,e_n\}$ be an orthonormal basis of $T_pM$ and $c=K(e_1,e_2)$. We will show that the sectional curvature of every $2$-plane in $T_pM$ is $c$. To this end, for each $i\in \{2,...,n\}$ denote by $V_i$ the linear span of $\{e_1,...,e_i\}$, so that

$$V_2\subset V_3\subset \cdots \subset V_{n-1}\subset V_n=T_pM.$$
Consider the set $\mathcal O\subset \{2, \dots, n\}$ defined by the property that $j\in \mathcal O$ if and only if
the sectional curvature    of every $2$-plane $\sigma \subset V_j$  satisfies $K(\sigma)=c$. Notice that $\mathcal O\neq \emptyset$, as $2\in \mathcal O$. Let now
 $\kappa=\max \mathcal O$. Clearly, Theorem \ref{point} ii) holds if
$\kappa=n$. Let's argue by contradiction and assume that
$2\leq \kappa \leq n-1.$ In particular,

\vskip8pt
\noindent ($\dagger$) \;\; $K(\sigma)=c$ \;\;$\forall$  $\sigma \in G_2( V_{\kappa})$, \;\; $\exists \; \tilde \sigma \in G_2(V_{\kappa+1}) : K(\tilde \sigma)\neq c$.
\vskip8pt
\noindent The possibilities for $\text{dim} (V_{\kappa}\cap \tilde\sigma)$ are, of course, $0,1$ or $2$.
The first case cannot occur otherwise  $\kappa+2=\text{dim} (V_{\kappa}\oplus \tilde\sigma)\leq \text{dim} V_{\kappa+1}=\kappa+1$. Since $ K(\tilde \sigma)\neq c$ one cannot have $\text{dim} (V_{\kappa}\cap \tilde\sigma)= 2$ either, and so the dimension of $V_{\kappa}\cap \tilde\sigma$ is $1$.

Let $w_1$ be a unit vector that generates  $V_{\kappa}\cap \tilde\sigma$, $w_2$ a unit vector such that
$\{w_1, w_2\}$ is an orthonormal basis of $\tilde \sigma$. Since $\text{dim} V_{\kappa}\geq 2$,  one can choose $w_3\in V_{\kappa}$ such that $\{w_1, w_3\}$ is an orthonormal set.
By ($\dagger$),
\begin{eqnarray} \label {c} K(w_1,w_3)=c \neq K(\tilde \sigma)=K(w_1, w_2). \end{eqnarray}
By Lemma \ref{onlyone} there exists a sequence $(q_k)$ in $M-\{p\}$ such that

$$\lim_{k \to \infty} q_k=p, \;\;\;
\lim_{k \to \infty}\frac{\nabla f(q_k)}{|\nabla  f(q_k)|} = w_1.$$
By  (\ref{conf}), for every convergent sequence of unit vectors $z_k\perp \nabla f(q_k)$, say $z_k\to z$, one has

\begin{eqnarray}  \label{near the end}
K(\frac{\nabla f(q_k)}{|\nabla f(q_k)|},z_k)=-\frac{\langle\nabla f,\nabla h\rangle}{|\nabla f|^2}(q_k).
\end{eqnarray}
The LHS of (\ref{near the end}) converges to $K(w_1, z)$. Manifestly, the RHS of (\ref{near the end})  is independent of $z_k$. It follows that there exists a constant $d$ such that
\vskip8pt
\noindent ($\dagger \dagger$) $K(w_1,z)=d \;\;\; \forall z\in w_1^{\perp}$, $|z|=1$.
\vskip8pt
\noindent  Taking  $z=w_2$ and $z=w_3$ in ($\dagger \dagger$),
\begin{eqnarray} \label{d} K(w_1,w_3)=d =K(w_1, w_2),\end{eqnarray}
a direct contradiction to (\ref{c}). This concludes the proof of Theorem \ref{point}.\qed

\vskip20pt

Xiaoyang Chen
\vskip3pt
School of Mathematical Sciences, Institute for Advanced Study
\vskip3pt
Tongji University
\vskip3pt
Shanghai, China
\vskip3pt
xychen100@tongji.edu.cn

\vskip20pt

 Francisco Fontenele
 \vskip3pt
Departamento de Geometria
\vskip3pt
Universidade Federal Fluminense
\vskip3pt
Niter$\acute{o}$i, RJ, Brazil
\vskip3pt
fontenele@mat.uff.br

\vskip20pt

  Frederico Xavier
  \vskip3pt
 Department of Mathematics
 \vskip3pt
 Texas Christian University
 \vskip3pt
 Fort Worth, TX, USA
 \vskip3pt
 f.j.xavier@tcu.edu

\end{document}